\documentclass[11pt]{article}
\usepackage[T1]{fontenc}
\usepackage[utf8]{inputenc}
\usepackage{lmodern,microtype}
\usepackage[a4paper,margin=28mm]{geometry}
\usepackage{amsmath,amssymb,amsthm,mathtools}
\usepackage{xcolor,enumitem}
\usepackage{authblk}

\usepackage[colorlinks=true,linkcolor=blue,citecolor=blue,urlcolor=blue,hypertexnames=false]{hyperref}
\hypersetup{pdftitle={One exponential or two? Invertibility tests for Jordan epimorphisms},pdfauthor={Catalin Badea}}
\allowdisplaybreaks
\numberwithin{equation}{section}
\setlist[enumerate]{itemsep=2pt,topsep=4pt}

\newtheorem{theorem}{Theorem}[section]
\newtheorem{proposition}[theorem]{Proposition}
\newtheorem{lemma}[theorem]{Lemma}
\newtheorem{corollary}[theorem]{Corollary}
\theoremstyle{definition}
\newtheorem{remark}[theorem]{Remark}

\newcommand{\C}{\mathbb C}
\newcommand{\Z}{\mathbb Z}
\newcommand{\T}{\mathbb T}
\newcommand{\ad}{\operatorname{ad}}
\newcommand{\one}{\mathbf 1}
\newcommand{\rad}{\operatorname{rad}}
\newcommand{\id}{\mathrm{id}}
\newcommand{\Bcal}{\mathcal B}
\newcommand{\spec}{\sigma}
\newcommand{\hsp}{\widehat\sigma}
\newcommand{\R}{\mathbb R}

\renewcommand{\Re}{\operatorname{Re}}
\renewcommand{\Im}{\operatorname{Im}}

\title{One exponential or two?\\ Invertibility tests for Jordan epimorphisms}
\author{Catalin Badea}
\affil{%
Department of Mathematics and Statistics, University of Reading, United Kingdom,\\
and UMR 8524 -- Laboratoire Paul Painlev\'e, Universit\'e de Lille, France\\
\texttt{c.badea@reading.ac.uk}, \texttt{cbadea@univ-lille.fr}\\
\url{https://www.reading.ac.uk/maths-and-stats/staff/catalin-badea}%
}
\date{}

\begin{document}
\maketitle

\begin{abstract} How many exponential factors must one test to recognise a Jordan
epimorphism? Let \(T\colon A\to B\) be a unital linear surjection
between complex Banach algebras, with \(B\) semisimple. Using a
spectral-isometry theorem of Aznay and Chabbabi, we show that it
suffices for \(T\) to send every product of two exponentials to an
invertible element to force $T$ to be Jordan. If only single exponentials are tested, \(T\) is a
Jordan epimorphism if and only if its kernel is a two-sided ideal, or
equivalently a left ideal. A resolvent estimate on a logarithmic spiral
shows that a single exponential suffices when the source is a
\(C^*\)-algebra or the group algebra \(\ell^1(G)\) of a discrete group.
A separate argument, which does not use the spectral-isometry theorem,
shows that a single exponential also suffices when the target is a
\(C^*\)-algebra with a faithful family of tracial states.
\end{abstract}



\section{Introduction}\label{sec:introduction}

\subsection{The Gleason--Kahane--\.Zelazko theorem and the
Kaplansky--Aupetit problem}
The Gleason--Kahane--\.Zelazko theorem characterises the multiplicative
linear functionals on a Banach algebra: a linear functional $\varphi$ on
a complex unital Banach algebra $A$ such that $\varphi(\one)=1$ and
$\varphi$ does not vanish on the invertible elements of $A$ is
multiplicative \cite{Gleason1967,KahaneZelazko1968,Zelazko1968}. The
proof gives more: it suffices that $\varphi$ does not vanish on
exponentials,
\begin{equation}\label{eq:scalar-exponential}
 \varphi(e^a)\ne0\quad(a\in A)
 \quad\Longrightarrow\quad
 \varphi(xy)=\varphi(x)\varphi(y)\quad(x,y\in A).
\end{equation}
For this formulation see, for instance, Jarosz \cite[Theorem~2]{Jarosz1997}
or \cite[Section~5.6]{BadKah}.

What replaces multiplicativity when the range is a noncommutative
algebra? Kaplansky \cite{Kaplansky1970} asked which unital linear maps
preserving invertibility are Jordan homomorphisms. In the form given
to the question by Aupetit \cite{AupetitPJM,Aupetit1991}, the
\emph{Kaplansky--Aupetit problem} reads: if $T\colon A\to B$ is a
unital linear surjection onto a semisimple Banach algebra $B$ which
sends \emph{all} invertible elements to invertible elements, is $T$ a
Jordan epimorphism, that is, does $T(x^2)=(Tx)^2$ hold for every
$x\in A$? Equivalently, $T$ should preserve the Jordan product
$(xy+yx)/2$. In a recent preprint, Aznay and Chabbabi
\cite{AznayChabbabi2026} have given a positive answer.\footnote{Theorem~\ref{thm:main},
the last assertion of Lemma~\ref{lem:E1-conjugation}, and all results
deduced from them depend on \cite{AznayChabbabi2026}; this concerns
Sections~\ref{sec:two}--\ref{sec:source} and
Theorem~\ref{thm:primitive-E1}. Proposition~\ref{prop:E1-hull},
Corollary~\ref{cor:E1-source-radical}, Section~\ref{sec:rank} and
Theorem~\ref{thm:tracial-E1} do not. As \cite{AznayChabbabi2026} is a
recent preprint, each use is indicated where it occurs.}
With \eqref{eq:scalar-exponential} in mind, we ask whether
invertible elements can be replaced by exponentials in the
Kaplansky--Aupetit problem, or at least by products of a fixed number
of exponentials.

\subsection{Notation and the two tests}
Throughout, $A$ and $B$ are complex, associative, unital Banach
algebras whose unit $\one\ne0$ has norm one. We write $A^{-1}$ for the
group of invertible elements, $\rad(A)$ for the Jacobson radical, and
$\spec_A(x)$ and $r_A(x)$ for the spectrum and the spectral radius of
$x\in A$. All linear maps are complex-linear; continuity is not
assumed.

Let $T\colon A\to B$ be a unital linear surjection, with $B$
semisimple. We compare two tests:
\begin{align}
 T(e^a)&\in B^{-1} &&(a\in A),\tag{E$_1$}\label{hyp:E1}\\
 T(e^ae^b)&\in B^{-1} &&(a,b\in A).\tag{E$_2$}\label{hyp:E2}
\end{align}
The two tests agree when the source is commutative, since then
$e^ae^b=e^{a+b}$. For a noncommutative source, a product of two
exponentials can lie at positive distance from the set of single
exponentials (Section~\ref{sec:rank}).

By \eqref{eq:scalar-exponential}, the single-exponential test
\eqref{hyp:E1} suffices when $B=\C$. Mathieu and Schulz observed that
their solution of the Kaplansky--Aupetit problem for $C^*$-algebra
targets with a faithful family of tracial states uses invertibility
only on products of two exponentials
\cite[Section~4, before Corollary~4.4]{MathieuSchulz2023}.
Theorem~\ref{thm:main} below shows that the two-exponential test
\eqref{hyp:E2} suffices for every semisimple target.

With a single exponential, the question becomes whether $\ker T$ is a
two-sided ideal (Theorem~\ref{thm:E1-ideal}). We establish this ideal
property by estimating the resolvent of $e^{zh}$ on a logarithmic
spiral, directly from the resolvent of a generator $h$. This shows that
\eqref{hyp:E1} suffices for every discrete group algebra $\ell^1(G)$
(Corollary~\ref{cor:E1-group-source}), even for $\ell^1(F_\infty)$, in
which the number of exponential factors needed to represent the
invertible elements of the identity component is unbounded. On the
target side, traces detect the defect $T(x^2)-(Tx)^2$. We end with a
description of what a counterexample would have to look like. Whether
one exponential suffices for arbitrary sources and semisimple targets
remains open.

\subsection{Principal results}\label{subsec:principal}
Throughout, \(T\colon A\to B\) is a unital linear surjection onto a
semisimple Banach algebra \(B\). Three results describe the general
picture.
\begin{itemize}
\item \emph{Two exponentials always suffice} (Theorem~\ref{thm:main}).
If \(T\) satisfies \eqref{hyp:E2}, then \(T\) is a bounded Jordan
epimorphism, its kernel is a closed two-sided ideal containing
\(\rad(A)\), and
\(\spec_B(Tx)=\spec_{A/\ker T}(x+\ker T)\subseteq\spec_A(x)\). 
\item \emph{One exponential and the kernel} (Theorem~\ref{thm:E1-ideal}).
If \(T\) satisfies \eqref{hyp:E1}, then \(T\) is a Jordan epimorphism if
and only if \(\ker T\) is a two-sided ideal. Since \(\ker T\) is
automatically a Lie ideal (Lemma~\ref{lem:E1-conjugation}), it is
equivalent to ask that \(\ker T\) be a left ideal. In particular, one
exponential suffices whenever \(T\) is injective.
\item \emph{One exponential for good sources}
(Theorems~\ref{thm:source-resolvent} and~\ref{thm:source-profile},
Corollaries~\ref{cor:source-cstar} and~\ref{cor:E1-group-source}). One
exponential suffices when \(A\) is generated by elements \(h\) with real
spectrum whose resolvents grow at most like \(|y|^{-1-\varepsilon}\),
for every \(\varepsilon>0\), at distance \(|y|\) from the real axis.
This covers unital \(C^*\)-algebras, algebras generated by
idempotents, all group algebras \(\ell^1(G)\) of discrete groups, and
weighted group algebras \(\ell^1(G,\omega)\) whose weight grows
subpolynomially along the powers of a generating set.
\end{itemize}
On the target side, Theorem~\ref{thm:tracial-E1} treats
\(C^*\)-algebras with a faithful family of tracial states, and
Theorem~\ref{thm:primitive-E1} and
Corollary~\ref{cor:primitive-obstruction} show that a counterexample
must produce a vector-valued trace onto an infinite-dimensional
noncommutative primitive algebra. Section~\ref{sec:rank} shows that
the number of exponential factors needed to \emph{test} a map and the
number needed to \emph{factorise} the elements of the identity
component of the invertible group are unrelated: for
\(\ell^1(F_\infty)\) the first is one and the second is unbounded.

\subsection{Organisation of the paper}
Section~\ref{sec:two} proves Theorem~\ref{thm:main},
Section~\ref{sec:kernel-core} the kernel criterion, and
Section~\ref{sec:source} the resolvent and group-algebra criteria.
Section~\ref{sec:rank} discusses exponential rank, and
Section~\ref{sec:targets} treats traces and the primitive obstruction.
A reader interested only in group algebras needs from
Section~\ref{sec:source} only Theorem~\ref{thm:source-profile},
Corollary~\ref{cor:source-growth} and Section~\ref{subsec:groups}.

\section{Why two exponentials suffice}\label{sec:two}

We first record the elementary fact that unital Jordan homomorphisms
preserve invertibility.

\begin{lemma}\label{lem:jordan-invertible}
Let $J$ be a unital Jordan homomorphism between unital algebras over
$\C$. If $x$ is invertible, then $J(x)$ is invertible and
$J(x^{-1})=J(x)^{-1}$.
\end{lemma}
\begin{proof}
Polarising $J(x^2)=J(x)^2$ shows that $J$ preserves $xy+yx$, and then
also $xyx$, because
$2xyx=(xy+yx)x+x(xy+yx)-(x^2y+yx^2)$. Put $u=J(x)$ and
$v=J(x^{-1})$. Applied to $x$ and $x^{-1}$, these identities give
$uvu=u$ and $uv+vu=2\cdot\one$. Hence $p=uv$ and $q=vu$ satisfy
$p^2=(uvu)v=p$, $q^2=v(uvu)=q$ and $p+q=2\cdot\one$. Substituting
$q=2\cdot\one-p$ in $q^2=q$ gives $p=\one$, and then $q=\one$.
\end{proof}

\begin{theorem}[The two-exponential test]\label{thm:main}
Let $T\colon A\to B$ be a unital linear surjection, where $B$ is
semisimple. If $T(e^ae^b)$ is invertible for all $a,b\in A$, then $T$
is a bounded Jordan epimorphism. Its kernel $I$ is a closed two-sided
ideal containing $\rad(A)$, and the induced map $A/I\to B$ is a
bicontinuous Jordan isomorphism. In particular,
\[
 \spec_B(Tx)=\spec_{A/I}(x+I)\subseteq\spec_A(x)\qquad(x\in A).
\]
\end{theorem}

We adapt the proof of Aznay and Chabbabi
\cite[Proposition~5.2]{AznayChabbabi2026}, which starts from the
inclusion $\spec_B(Tx)\subseteq\spec_A(x)$. The point is that their
normalised maps require only a bound on the spectral radius, and that
products of two exponentials suffice to obtain it.

\begin{proof}
\emph{Step~1: boundedness and the radical.} If $|\lambda|>r_A(x)$,
then $\lambda\one-x=e^c$ with
\begin{equation}\label{eq:log-series}
 c=(\log\lambda)\one-\sum_{n\ge1}\frac{(x/\lambda)^n}{n},
\end{equation}
for any scalar logarithm of $\lambda$; the series converges absolutely
by the spectral radius formula. Thus $\lambda\one-Tx=T(e^ce^0)$ is
invertible, and $r_B(Tx)\le r_A(x)$. By the automatic continuity
theorem for spectrally bounded surjections onto semisimple Banach
algebras \cite{Aupetit1982}, \cite[Theorem~2.1]{Mathieu2009}, $T$ is
bounded. Such a surjection also maps the radical into the radical
\cite[Proposition~2.11]{MathieuSchick2002}; hence $\rad(A)\subseteq I$.

\emph{Step~2: normalised maps.} Fix $a\in A$. For $z\in\C$ put
$u_z=T(e^{za})$, which is invertible, and define
\[
 \Phi_z(x)=u_z^{-1}T(e^{za}x),\qquad
 \Psi_z(x)=T(xe^{za})u_z^{-1}\qquad(x\in A).
\]
Both maps are unital linear surjections, and they depend
holomorphically on $z$. With $c$ as in \eqref{eq:log-series},
\[
 \lambda\one-\Phi_z(x)=u_z^{-1}T(e^{za}e^c),\qquad
 \lambda\one-\Psi_z(x)=T(e^ce^{za})u_z^{-1},
\]
and both elements are invertible. Hence
$r_B(\Phi_z(x))\le r_A(x)$ and $r_B(\Psi_z(x))\le r_A(x)$ for every
$z\in\C$. By Vesentini's theorem \cite{Vesentini1968}, these spectral
radii are subharmonic functions of $z$; being bounded above on $\C$,
they are constant by Liouville's theorem for subharmonic functions.
Since $\Phi_0=\Psi_0=T$,
\begin{equation}\label{eq:constant-radius}
 r_B(\Phi_z(y))=r_B(\Psi_z(y))=r_B(Ty)\qquad(y\in A,\ z\in\C).
\end{equation}

\emph{Step~3: the kernel is an ideal.} We use Zem\'anek's
characterisation of the radical \cite[Theorem~5.3.1]{Aupetit1991}:
\begin{equation}\label{eq:radical-criterion}
 q\in\rad(B)\quad\Longleftrightarrow\quad
 r_B(y+q)=r_B(y)\quad(y\in B).
\end{equation}
Fix $j\in I$. By \eqref{eq:constant-radius},
$r_B(\Phi_z(x)+\Phi_z(j))=r_B(T(x+j))=r_B(Tx)=r_B(\Phi_z(x))$ for every
$x\in A$, and $\Phi_z(x)$ ranges over $B$ as $x$ ranges over $A$. By
\eqref{eq:radical-criterion} and semisimplicity, $\Phi_z(j)=0$, and in
the same way $\Psi_z(j)=0$. Thus $T(e^{za}j)=T(je^{za})=0$ for all
$z$, and differentiation at $z=0$ gives $aj,ja\in I$. Therefore $I$ is
a closed two-sided ideal. Every product of two exponentials in $A/I$
lifts to a product of two exponentials in $A$, so the induced map
$A/I\to B$ satisfies the same hypothesis; it is a bounded linear
bijection, with bounded inverse by the open mapping theorem. Replacing
$A$ by $A/I$, we may therefore assume that $T$ is bijective, with
bounded inverse. The identities \eqref{eq:constant-radius} remain
valid, since the normalised maps of the induced map are induced by
$\Phi_z$ and $\Psi_z$.

\emph{Step~4: the transported product.} Transport the multiplication
of $A$ to $B$ by setting $b\star c=T(T^{-1}(b)T^{-1}(c))$. By
\eqref{eq:constant-radius}, the maps $V_z=\Phi_zT^{-1}$ are unital
bijective spectral isometries of $B$; by
\cite[Theorem~1.1]{AznayChabbabi2026}, they are Jordan isomorphisms.
Put $b=Ta$. Since $V_0=\id$ and $\frac{d}{dz}u_z^{-1}\big|_{z=0}=-b$,
\[
 \frac{d}{dz}V_z(c)\Big|_{z=0}=-bc+T\bigl(aT^{-1}(c)\bigr)
                              =b\star c-bc=:D_b(c).
\]
Differentiating $V_z(c^2)=V_z(c)^2$ at $z=0$ shows that $D_b$ is a
Jordan derivation of $B$. As $B$ is semiprime, $D_b$ is a derivation
\cite{Bresar1988}. This holds for every $b\in B$, since $a$ was
arbitrary.

\emph{Step~5: the Jordan identity.} This algebraic step uses only the
left normalisations. Fix $b\in B$, write $L_c(x)=cx$ and put $D=D_b$, so
that $c\mapsto b\star c$ is $L_b+D$. Associativity of $\star$ gives
$(L_b+D)^2=L_{b\star b}+D_{b\star b}$. Substituting
$DL_b=L_bD+L_{D(b)}$ and $b\star b=b^2+D(b)$, we obtain
$D_{b\star b}=2L_bD+D^2$. Since $D$ is a derivation, the Leibniz
defects of $L_bD$ and of $D^2$ are $(x,y)\mapsto[b,x]D(y)$ and
$(x,y)\mapsto2D(x)D(y)$. As $D_{b\star b}$ is a derivation, it follows
that
\[
 \bigl([b,x]+D(x)\bigr)D(y)=0\qquad(x,y\in B).
\]
Put $x=b$ and $q=D(b)$. Then $qD(y)=0$ for every $y$, and replacing $y$
by $yb$ gives $qyq=0$ for every $y\in B$. Thus $BqB$ is a two-sided
ideal with square zero, and $q=0$ by semiprimeness. Hence
$b\star b=b^2$ for every $b\in B$, which is the Jordan identity for $T$.

Finally, $T$ and $T^{-1}$ are unital Jordan isomorphisms, so both
preserve invertibility by Lemma~\ref{lem:jordan-invertible}. This gives
the spectral identity; the inclusion holds because the quotient map
$A\to A/I$ is a unital homomorphism.
\end{proof}

\section{One exponential and the kernel}\label{sec:kernel-core}

The normalised maps of Section~\ref{sec:two} need products of two
exponentials. Under \eqref{hyp:E1} alone, conjugation by exponentials
still controls commutators. We first record what a single exponential
detects in the spectrum.

\begin{proposition}\label{prop:E1-hull}
Let $T\colon A\to B$ be unital and linear and satisfy \eqref{hyp:E1}.
Writing $\hsp(x)$ for the polynomially convex hull of the spectrum,
one has
\[
 \hsp_B(Tx)\subseteq\hsp_A(x),\qquad r_B(Tx)\le r_A(x)
 \qquad(x\in A).
\]
If $T$ is onto and $B$ semisimple, it is bounded.
\end{proposition}
\begin{proof}
Let $\lambda\notin\hsp_A(x)$. Then $\lambda$ lies in the unbounded
component of $\C\setminus\spec_A(x)$, so there is a path in this
component from $\lambda$ to some $\mu$ with $|\mu|>\|x\|$. Let
$C=\{x\}''$ be the bicommutant of $x$; it is a closed commutative
unital subalgebra of $A$, and it contains the inverse of every element
of $C$ which is invertible in $A$. Hence the elements $\xi\one-x$ along
the path are invertible in $C$, and the path joins $\lambda\one-x$ to
$\mu\one-x$ in $C^{-1}$. By \eqref{eq:log-series}, $\mu\one-x$ has a
logarithm in $C$. In a commutative unital Banach algebra the identity
component of the invertible group consists of single exponentials:
subdivide a path into increments close to $\one$, take their
logarithms and add them. Thus $\lambda\one-x=e^c$ for some $c\in C$,
and $\lambda\one-Tx=T(e^c)$ is invertible. This proves
$\spec_B(Tx)\subseteq\hsp_A(x)$; since $\hsp_A(x)$ is polynomially
convex, the hull inclusion and the radius bound follow. The automatic
continuity theorem used in Step~1 of the proof of
Theorem~\ref{thm:main} gives the last assertion.
\end{proof}

\begin{corollary}\label{cor:E1-source-radical}
Let $T\colon A\to B$ be a unital linear surjection onto a semisimple
Banach algebra satisfying \eqref{hyp:E1}. Then $T$ annihilates
$\rad(A)$. The induced map $A/\rad(A)\to B$ satisfies \eqref{hyp:E1},
and it is a Jordan epimorphism if and only if $T$ is one.
\end{corollary}
\begin{proof}
The first assertion follows from the radius bound of
Proposition~\ref{prop:E1-hull} and
\cite[Proposition~2.11]{MathieuSchick2002}. Exponentials in
$A/\rad(A)$ lift to exponentials in $A$, and the square identity passes
to and from the quotient.
\end{proof}

\begin{lemma}[Conjugation on the range]\label{lem:E1-conjugation}
Let $T\colon A\to B$ be a unital linear surjection onto a semisimple
Banach algebra, and suppose that there is $M\ge0$ with
$r_B(Tx)\le M r_A(x)$ for every $x\in A$. Then $\ker T$ is a closed
Lie ideal: $[A,\ker T]\subseteq\ker T$. Moreover, for every $a\in A$
the formula $D_a(Tx)=T([a,x])$ defines a bounded derivation of $B$.
\end{lemma}
\begin{proof}
The map $T$ is bounded by automatic continuity, so $\ker T$ is closed.
Fix $a\in A$ and put $C_z(x)=T(e^{za}xe^{-za})$. Since the spectral
radius is invariant under similarity,
$r_B(C_z(x))\le Mr_A(x)$ for every $z\in\C$; by the theorems of
Vesentini and Liouville, $r_B(C_z(x))=r_B(C_0(x))=r_B(Tx)$. For
$j\in\ker T$, comparing $x+j$ with $x$ gives
$r_B(C_z(x)+C_z(j))=r_B(C_z(x))$ for every $x\in A$. As $C_z$ is onto,
\eqref{eq:radical-criterion} gives $C_z(j)=0$, and differentiation at
$z=0$ gives $T([a,j])=0$. This proves the first assertion.

Consequently $U_z(Tx)=C_z(x)$ well defines a unital bijective spectral
isometry $U_z$ of $B$, which is a Jordan isomorphism by
\cite[Theorem~1.1]{AznayChabbabi2026}. Since $U_0=\id$ and
$\frac{d}{dz}U_z\big|_{z=0}=D_a$, differentiating
$U_z(c^2)=U_z(c)^2$ at $z=0$ shows that $D_a$ is a Jordan derivation,
and hence a derivation \cite{Bresar1988}. It is bounded, since it is
conjugate, through the induced bijection $A/\ker T\to B$ (bounded with
bounded inverse), to the operator induced by $\ad_a$ on $A/\ker T$.
\end{proof}

\begin{theorem}[The ideal-kernel criterion]\label{thm:E1-ideal}
Let $T\colon A\to B$ be a unital complex-linear surjection onto a
semisimple Banach algebra satisfying \eqref{hyp:E1}. Then $T$ is a
Jordan epimorphism if and only if $\ker T$ is a two-sided ideal.
In particular, the single-exponential test suffices for a bijection.
\end{theorem}
\begin{proof}
A unital Jordan epimorphism satisfies \eqref{hyp:E2} by
Lemma~\ref{lem:jordan-invertible}, so Theorem~\ref{thm:main} gives
the ideal property. Conversely, suppose $I=\ker T$ is a two-sided
ideal. It is closed by Proposition~\ref{prop:E1-hull}. The induced
map $A/I\to B$ satisfies \eqref{hyp:E1}, since exponentials lift;
we may thus suppose that $T$ is bijective.

Fix $x\in A^{-1}$. Following Costara and Repov\v{s}
\cite[proof of Theorem~2]{CostaraRepovs2010}, consider the
bicommutant $C=\{x\}''$. It is the centre of the commutant
$\{x\}'$, hence a closed commutative unital algebra containing
$x$ and $x^{-1}$. By Proposition~\ref{prop:E1-hull} and
Lemma~\ref{lem:E1-conjugation}, the operators $D_a=T\ad_a T^{-1}$ are
bounded derivations of $B$, and, as $T$ is injective,
\[
 E=T(C)=\bigcap_{a\in\{x\}'}\ker D_a .
\]
Thus $E$ is a closed unital subalgebra of $B$. It is inverse-closed
in $B$, since $D_a(v^{-1})=-v^{-1}D_a(v)v^{-1}$ for every $v\in B^{-1}$.
Consequently $T(e^c)\in E^{-1}$ for all $c\in C$, by \eqref{hyp:E1}.

Let $q\colon E\to E/\rad(E)$ be the quotient map. The unital
surjection $qT|_C\colon C\to E/\rad(E)$ satisfies \eqref{hyp:E2},
because $C$ is commutative and $e^ce^d=e^{c+d}$. By
Theorem~\ref{thm:main} it is a Jordan epimorphism, and by
Lemma~\ref{lem:jordan-invertible} the element $q(Tx)$ is invertible in
$E/\rad(E)$. Invertibility lifts modulo the radical, so
$Tx\in E^{-1}\subseteq B^{-1}$. We have shown that $T$ maps every
invertible element to an invertible element; in particular it
satisfies \eqref{hyp:E2}, and Theorem~\ref{thm:main} completes the
proof.
\end{proof}

The hypothesis of Theorem~\ref{thm:E1-ideal} can be weakened, because the
Lie ideal property is automatic. Under \eqref{hyp:E1} the bound
$r_B(Tx)\le r_A(x)$ holds by Proposition~\ref{prop:E1-hull}, so
Lemma~\ref{lem:E1-conjugation} applies and $[A,\ker T]\subseteq\ker T$.
Since $ja=aj-[a,j]$, the kernel is a two-sided ideal as soon as it is a left
ideal. Thus, under \eqref{hyp:E1},
\[
 T\ \text{is Jordan}
 \quad\Longleftrightarrow\quad
 \ker T\ \text{is a left ideal}
 \quad\Longleftrightarrow\quad
 T(aj)=0\ \ (a\in A,\ j\in\ker T).
\]
The resolvent criterion of Section~\ref{sec:source} establishes this
left ideal property and then invokes Theorem~\ref{thm:E1-ideal}.

The criterion concerns the Jordan identity, not equality of spectra in
the original source. For example, the projection $\C^2\to\C$ onto the
first coordinate is a unital homomorphism and satisfies both tests, but
it maps $(0,1)$, whose spectral radius is one, to zero.

Surjectivity cannot be omitted, even for the full invertibility test.
The unital map
\[
 (s,t)\longmapsto
 \begin{pmatrix}s&t-s\\0&s\end{pmatrix},\qquad \C^2\longrightarrow M_2(\C),
\]
maps invertible elements to invertible elements, but it maps the
idempotent $(0,1)$ to a nonzero matrix with square zero.

\begin{remark}[Testing modulo an ideal]
Both Theorem~\ref{thm:main} and Theorem~\ref{thm:E1-ideal} require a
semisimple target, but invertibility may be tested modulo an ideal.
Let $J$ be a proper closed two-sided ideal of $B$ with $B/J$
semisimple, let $\pi\colon B\to B/J$ be the quotient map, and let
$T\colon A\to B$ be a linear map, not necessarily unital, such that
$\pi T$ is onto and $\pi T(e^a)$, respectively $\pi T(e^ae^b)$, is
invertible in $B/J$ for all $a,b\in A$. Taking $a=b=0$ shows that
$u=\pi(T\one)$ is invertible, and $S=u^{-1}\pi T$ is a unital surjection
onto $B/J$ which satisfies the same test, with $\ker S=T^{-1}(J)$. Under
the two-exponential test, $S$ is a Jordan epimorphism by
Theorem~\ref{thm:main}; under the single-exponential test, it is one if
and only if $T^{-1}(J)$ is a two-sided ideal, by
Theorem~\ref{thm:E1-ideal}. For instance, if $H$ is an
infinite-dimensional Hilbert space, $B=\Bcal(H)$ and $J=\mathcal K(H)$,
then $B/J$ is a $C^*$-algebra, hence semisimple, and by Atkinson's
theorem the tests become Fredholm tests: $T(e^a)$, respectively
$T(e^ae^b)$, is a Fredholm operator.
\end{remark}

\section{A source criterion from the resolvent}\label{sec:source}

For $h\in A$ write $R_h(\xi)=(\xi\one-h)^{-1}$, defined on the resolvent set of
$h$. The next criterion concerns the behaviour of $R_h$ near the real
axis. The constants may depend on the generator.

\begin{theorem}[A resolvent test for generators]\label{thm:source-resolvent}
Suppose that $A$ is generated as a closed unital algebra by elements
$h$ with real spectrum such that, for each $\varepsilon>0$,
\begin{equation}\label{eq:almost-linear-resolvent}
 \|R_h(x+iy)\|\le C_{h,\varepsilon}|y|^{-1-\varepsilon}
 \qquad (x\in\R,\quad 0<|y|\le1).
\end{equation}
If $B$ is semisimple, every unital linear surjection $T\colon A\to B$
satisfying \eqref{hyp:E1} is a bounded Jordan epimorphism.
In particular, a bound $\|R_h(x+iy)\|\le C_h/|y|$ suffices.
\end{theorem}

\begin{corollary}\label{cor:source-cstar}
A single exponential suffices if the source is a unital $C^*$-algebra
or a Banach algebra generated by idempotents. The latter class includes
every Banach algebra containing a unital copy of $M_n(\C)$, $n\ge2$.
\end{corollary}
\begin{proof}
For a self-adjoint element of a $C^*$-algebra, the continuous functional
calculus gives $\|R_h(x+iy)\|\le1/|y|$, and self-adjoint elements span
the algebra. For an idempotent $p$ use
$R_p(\xi)=\xi^{-1}(\one-p)+(\xi-1)^{-1}p$; since $|y|\le|\xi|$ and
$|y|\le|\xi-1|$ for $\xi=x+iy$, this gives
$|y|\,\|R_p(x+iy)\|\le\|\one-p\|+\|p\|$. The spectrum of $p$ is
contained in $\{0,1\}$, and both assertions follow from
Theorem~\ref{thm:source-resolvent}.

For the last assertion, take matrix units $(e_{ij})$ with
$\sum_i e_{ii}=\one$. Each off-diagonal corner $x=e_{ii}xe_{jj}$
is a difference of idempotents: $x=(e_{ii}+x)-e_{ii}$.
Each diagonal corner $d=e_{ii}de_{ii}$ is the product
$(de_{ij})e_{ji}$ of two off-diagonal elements for any $j\ne i$.
Since every element is the sum of its corners, the idempotents
generate the algebra even algebraically.
\end{proof}

The bound $C_h/|y|$ also holds for a generator $h$ with a bounded
absolutely continuous functional calculus on a compact real interval
$I$, that is, a bounded unital homomorphism $AC(I)\to A$ which maps the
coordinate function to $h$; here $\|f\|_{AC}=\|f\|_\infty+\int_I |f'|$.
Indeed, $f(s)=(x+iy-s)^{-1}$ has supremum norm at most $1/|y|$ and
\[
 \int_I|f'(s)|\,ds\le \int_{\R}\frac{ds}{(s-x)^2+y^2}
 =\frac{\pi}{|y|}.
\]
The image of $f$ is the resolvent of $h$ at $x+iy$; the same argument
for $\xi\notin I$ shows that $\spec_A(h)\subseteq I$. This includes the
well-bounded operators discussed by Barnes \cite[Theorem~7(C)]{Barnes}.
Here it matters that the hypothesis concerns the resolvent rather than
the imaginary exponentials: for the coordinate function $h$ in
$AC([0,1])$ one has $\|e^{ith}\|_{AC}=1+|t|$, so the growth condition of
Corollary~\ref{cor:source-growth} below fails, whereas
Theorem~\ref{thm:source-resolvent} applies.

\subsection{The contour estimate}
We prove Theorem~\ref{thm:source-resolvent} in a form that allows the
resolvent to grow at different rates on different scales. For a generator
with real spectrum put
\[
 D_h(L)=\max\left\{1,\sup_{\substack{x\in\R\\ L^{-1}\le |y|\le1}}
          |y|\,\|R_h(x+iy)\|\right\}\qquad(L\ge1).
\]
This is finite for each $L$: on the region $L^{-1}\le|y|\le1$ the point
$x+iy$ is at distance at least $L^{-1}$ from the spectrum, the
resolvent is continuous there, and $\|R_h(\xi)\|\le2|\xi|^{-1}$ for
$|\xi|\ge2\|h\|$ by the Neumann series. Clearly $D_h$ is
nondecreasing.

\begin{theorem}\label{thm:source-profile}
The conclusion of Theorem~\ref{thm:source-resolvent} holds if its
generators have real spectrum and each satisfies
\begin{equation}\label{eq:resolvent-scales}
 D_h(L)\le C_h L^{p_h}\quad(L\ge1)
 \quad\text{for some }p_h\ge0,\qquad
 \liminf_{L\to\infty}\frac{\log D_h(L)}{\log L}=0 .
\end{equation}
\end{theorem}
Indeed \eqref{eq:almost-linear-resolvent} gives
$D_h(L)\le\max\{1,C_{h,\varepsilon}\}L^\varepsilon$ for every
$\varepsilon>0$, so the theorem implies Theorem~\ref{thm:source-resolvent}.

\begin{proof}
Fix one generator $h$, and for $z\in\C$ put
\[
 U_z=e^{zh},\qquad \Phi_z(x)=T(U_z)^{-1}T(U_zx)\quad(x\in A).
\]
By Proposition~\ref{prop:E1-hull}, $T$ is bounded and $T(U_z)$ is
invertible, so the maps $\Phi_z$ are unital surjections, entire in $z$
for the operator norm. Our aim is to prove
that $r_B(\Phi_z(x))$ is constant in $z$. For this we need a region of the
$\lambda$-plane on which $U_z(\lambda\one-x)$ is a single exponential,
uniformly in $\lambda$ there.

Fix $z$ with $\Re z\ne0$ and abbreviate $U=U_z$. Write $r=|z|$,
$k=\Im z/\Re z$ and $\kappa=r/|\Re z|$, and let
$\Gamma_k=\{-e^{zs_0}:s_0\in\R\}$.
This spiral is disjoint from $\spec_A(U)\subseteq e^{z\R}$: if
$-e^{zs_0}=e^{zs}$ with $s,s_0\in\R$, comparing moduli gives $s=s_0$,
which is absurd.
We claim that uniformly for $w\in\Gamma_k$,
\begin{equation}\label{eq:direct-spiral}
 \|(w\one-U)^{-1}U\|
 \le C_h\kappa^2[1+\log(1+r)]D_h(10\max\{1,\kappa r\}).
\end{equation}
Write $w=-e^{zs_0}$. The expression on the left is the norm of
$F(h)$ for $F(\xi)=(1+e^{z(\xi-s_0)})^{-1}-1$. Choose $b>0$
with $\spec_A(h)\subset(-b,b)$. Use a Cauchy contour consisting of
the graphs $x\pm i\eta(x)$ on $[-b,b]$, joined by vertical segments,
where
\[
 \eta(x)=\tfrac1{10}\min\{1,(|x-s_0|+r^{-1})/\kappa\}.
\]
This contour narrows near $s_0$, where a fixed strip might contain
poles of $F$. To see that $|F|\le C\kappa$ throughout the contour
and its interior, put $d=x-s_0$ and $v=z(d+iy)$, with $|y|\le\eta(x)$.
Note that $\kappa\ge1$ and $|F|=|e^v/(1+e^v)|$. If $r|d|\le1$, then
$r|y|\le r\eta(x)\le(r|d|+1)/10\le2/10$, so $|v|\le1+2/10<\pi/2$; hence
$\Re(e^v)>0$, $|1+e^v|>1$ and $|F|\le e^{6/5}$. Otherwise
$r^{-1}<|d|$, so $|y|\le|d|/(5\kappa)$ and
$|\Re v|\ge r|d|/\kappa-r|y|\ge c:=(4/5)r|d|/\kappa\ge4/(5\kappa)$. Since
$|1+e^v|\ge|1-e^{\Re v}|$, we obtain
$|F|\le\max\{e^c/(e^c-1),\,e^{-c}/(1-e^{-c})\}\le1+1/c\le 3\kappa$.
In particular $1+e^v$ does not vanish, and no pole of $F$ is
enclosed.

The integrals over the vertical ends are bounded by $C_h\kappa$,
since the ends lie in a fixed compact subset of the resolvent set of
$h$.
On either graph, $|d\xi|\le2\,dx$ and
$\eta(x)^{-1}\le10\max\{1,\kappa r\}$, so that
$\|R_h(x\pm i\eta(x))\|\le D_h(10\max\{1,\kappa r\})/\eta(x)$. Moreover,
\[
 \int_{-b}^{b}\frac{dx}{\eta(x)}
 \le C_b[1+\kappa\log(1+2br)],
\]
since the integral of $(|x-s_0|+r^{-1})^{-1}$ over this interval
is at most $2\log(1+2br)$. Thus the Cauchy integral for $F(h)$
gives \eqref{eq:direct-spiral}.

Write $p=p_h$ and take $L_j\to\infty$ with $\log D_h(L_j)=o(\log L_j)$.
We estimate \eqref{eq:direct-spiral} for those $z$ with $\Re z\ne0$ and
$|z|\le L_j^{1/4}$, so that $r\le L_j^{1/4}$. If $\kappa\le L_j^{1/4}$,
then the argument of $D_h$ in \eqref{eq:direct-spiral} is
$10\max\{1,\kappa r\}\le10L_j^{1/2}\le L_j$ for large $j$, and we use
the monotonicity of $D_h$. If instead
$\kappa>L_j^{1/4}$, then $r\le L_j^{1/4}<\kappa$, so $\kappa r\le\kappa^2$
and the polynomial bound in \eqref{eq:resolvent-scales} gives
$D_h(10\max\{1,\kappa r\})\le C'_h\kappa^{2p}$. In both cases the left-hand
side of \eqref{eq:direct-spiral} is bounded by
\[
 M_j(z):= C''_h D_h(L_j)\kappa^{2+2p}[1+\log(1+L_j)]\ge1 .
\]

For $|\lambda|>M_j(z)\|x\|$, the element $Y=U(\one-x/\lambda)$
is invertible. Its spectrum misses $\Gamma_k$, since
\[
 w\one-Y=(w\one-U)
     \bigl(\one+(w\one-U)^{-1}Ux/\lambda\bigr)\in A^{-1}
 \qquad(w\in\Gamma_k).
\]
On the complement of $\Gamma_k\cup\{0\}$ there is a holomorphic
logarithm: choose the argument of $w$ between $k\log|w|-\pi$ and
$k\log|w|+\pi$. Thus $Y$ has a logarithm, and adding a scalar
logarithm of $\lambda$ gives one for $\lambda Y$. By \eqref{hyp:E1},
$\lambda\one-\Phi_z(x)=T(U)^{-1}T(\lambda Y)$ is invertible.
We have proved that $r_B(\Phi_z(x))\le M_j(z)\|x\|$ whenever
$\Re z\ne0$ and $|z|\le L_j^{1/4}$. The imaginary axis, where this bound
is not available, meets each circle in at most two points, and so does
not affect the circle means below.

For $x\ne0$ the function $u(z)=\log r_B(\Phi_z(x))$ is subharmonic
\cite{Vesentini1968}, with $u\equiv-\infty$ allowed. For
$|z_0|\le L_j^{1/4}/4$ take the circle about $z_0$ of radius
$\rho=L_j^{1/4}/2$. It lies in the disc of the preceding estimate,
and the mean of $\log\kappa$ on it is $\log2$. Indeed, Jensen's
formula on $|\zeta|=1$ gives the means of $\log|z_0+\rho\zeta|$
and $\log|\Re(z_0+\rho\zeta)|$ as $\log\rho$ and $\log(\rho/2)$; for
the latter, write
$2\Re(z_0+\rho\zeta)=\rho\bar\zeta(\zeta^2+2(\Re z_0/\rho)\zeta+1)$ on
the circle, and note that the two roots of this quadratic lie on the
circle. The logarithmic singularities are integrable. The submean
inequality for $u$, together with the bound
$u(z)\le\log M_j(z)+\log\|x\|$, yields
\[
 \sup_{|z_0|\le L_j^{1/4}/4}u(z_0)
 \le \log D_h(L_j)+O_{h,x}(1+\log\log L_j)=o(\log L_j).
\]
Unless $u\equiv-\infty$, the maximum of $u$ on the disc of radius $R$
is a nondecreasing convex function of $\log R$. If it were unbounded,
it would grow at least like $\sigma\log R$ for some $\sigma>0$, which
contradicts this estimate along $L_j\to\infty$. Thus $u$ is bounded
above on $\C$, and hence constant. Evaluation at zero gives
\begin{equation}\label{eq:source-constant-radius}
              r_B(\Phi_z(x))=r_B(Tx)\qquad(x\in A,\ z\in\C).
\end{equation}

For $j_0\in\ker T$, apply this equality to $x+j_0$ and to $x$.
Surjectivity of $\Phi_z$ gives $r_B(y+\Phi_z(j_0))=r_B(y)$ for all
$y\in B$, so $\Phi_z(j_0)=0$ by \eqref{eq:radical-criterion} and
semisimplicity. Differentiating $T(e^{zh}j_0)=0$ at $z=0$ gives
$hj_0\in\ker T$. The set of $a\in A$ with $a\ker T\subseteq\ker T$ is
a closed unital subalgebra of $A$, because $\ker T$ is closed; it
contains every generator, and hence equals $A$. Thus $\ker T$ is a left
ideal. It is a Lie ideal by Lemma~\ref{lem:E1-conjugation}, and
therefore a two-sided ideal. Theorem~\ref{thm:E1-ideal} completes the
proof.
\end{proof}

\begin{remark}[A uniform local test]\label{rem:uniform-local}
There is also a local form of the argument. Suppose that \(T\) satisfies
\eqref{hyp:E1} and that \(S\) generates \(A\) as a closed unital algebra.
It suffices that for each \(h\in S\) there be an \(\varepsilon_h>0\)
such that \(T(e^{zh}e^b)\) is invertible for every \(z\in\C\) and every
\(\|b\|<\varepsilon_h\). Indeed, if
\(|\lambda|>\|x\|/(1-e^{-\varepsilon_h})\) and \(\lambda\ne0\), the
logarithm series gives \(\|\log(\one-x/\lambda)\|<\varepsilon_h\).
Thus
\(\lambda\one-\Phi_z(x)=\lambda T(e^{zh})^{-1}T(e^{zh}e^{b})\), with
\(b=\log(\one-x/\lambda)\), is invertible for all \(z\), and
\(r_B(\Phi_z(x))\) is bounded on the plane. Subharmonicity gives
\eqref{eq:source-constant-radius}, and the kernel argument above applies.
Conversely, every unital Jordan epimorphism passes this test, by
Lemma~\ref{lem:jordan-invertible}. Note that the neighbourhood of zero
in the hypothesis must not depend on \(z\).
\end{remark}

\begin{corollary}[Growth of the imaginary exponential groups]\label{cor:source-growth}
The conclusion of Theorem~\ref{thm:source-resolvent} holds if $A$ is
generated by elements $h$ for which
$\Omega_h(R)=\max_{|t|\le R}\|e^{ith}\|$ has polynomial growth and
\[
 \liminf_{R\to\infty}\frac{\log\Omega_h(R)}{\log R}=0.
\]
In particular it suffices that $\|e^{ith}\|=O_{h,\varepsilon}
((1+|t|)^\varepsilon)$ for every $\varepsilon>0$. The exponents,
constants and sequences of scales need not be the same for different
generators.
\end{corollary}
\begin{proof}
Suppose that $\Omega_h(R)\le C(1+R)^p$. Then
$\spec_A(h)\subset\R$: if $\mu\in\spec_A(h)$, the spectral mapping
theorem gives $e^{-t\Im\mu}\le r_A(e^{ith})\le\|e^{ith}\|$ for all
$t\in\R$, which forces $\Im\mu=0$. For $\nu=|y|>0$
the Laplace formula \cite[Theorem~1]{Barnes} gives
\[
 \nu\|R_h(x+iy)\|\le
 \nu\int_0^\infty e^{-\nu t}\|e^{-i\,\mathrm{sgn}(y)th}\|\,dt
\]
uniformly in $x$. (For Barnes's operator formulation one may use
left multiplication by $h$ on $A$.) For $0<\nu\le1$ this is
$O(\nu^{-p})$, so $D_h(L)=O(L^p)$. Choose $R_j\to\infty$ with
$\log\Omega_h(R_j)=o(\log R_j)$. For $R_j^{-1/2}\le\nu\le1$ the
integral including its prefactor is at most $\Omega_h(R_j)$ up to
$t=R_j$. Its tail is $O(e^{-\sqrt{R_j}/2}\nu^{-p})=o(1)$, by
retaining half of the exponential decay under the integral. Hence
$D_h(\sqrt{R_j})\le\Omega_h(R_j)+o(1)$. Apply
Theorem~\ref{thm:source-profile}.
\end{proof}

\subsection{Discrete group algebras}\label{subsec:groups}

A point mass of infinite order need not have a logarithm in its group
algebra. It can, however, be recovered from two real-valued functions
on the circle whose imaginary exponentials grow slowly in the Wiener norm.

\begin{corollary}[Discrete group algebras]\label{cor:E1-group-source}
Let $G$ be any discrete group. Every unital complex-linear surjection
$T:\ell^1(G)\to B$ onto a semisimple Banach algebra satisfying
\eqref{hyp:E1} is a bounded Jordan epimorphism. The same holds for
$\ell^1(G,\omega)$ when $\omega\ge1$ is a submultiplicative weight
with $\omega(e)=1$ and $G$ has a generating set $S$ such that
\[
 \omega(g^n)=O_{g,\varepsilon}\bigl((1+|n|)^\varepsilon\bigr)
 \qquad(g\in S,\ n\in\Z,\ \varepsilon>0).
\]
\end{corollary}
\begin{proof}
Identify $\T$ with $\R/2\pi\Z$. Let $b(x)=|x|$ on $[-\pi,\pi]$,
extended periodically, and put $c(x)=b(x-\pi/2)$. Then
\begin{equation}\label{eq:circle-cosines}
                    e^{ix}=\cos b(x)+i\cos c(x).
\end{equation}
For $0\le p<1$, write $W_p$ for the weighted Wiener algebra with norm
$\|f\|_{W_p}=\sum_{n\in\Z}(1+|n|)^p|\widehat f(n)|$.
Both $b$ and $c$ belong to $W_p$, since their Fourier coefficients
are $O((1+|n|)^{-2})$. We shall use the estimate
\[
 \|e^{itb}\|_{W_p}=\|e^{itc}\|_{W_p}
             \le C_p(1+|t|)^p\log(2+|t|)\qquad(t\in\R).
\]
For $p=0$ this is Kahane's logarithmic estimate for piecewise linear
phases \cite[Ch.~VI, \S2]{Kahane1970}. For the weighted estimate, write
the $n$th Fourier coefficient as
\[
 a_n(t)=\frac1{2\pi}\int_0^\pi
              \bigl(e^{i(t-n)x}+e^{i(t+n)x}\bigr)\,dx .
\]
For $|n|\le2|t|+2$, bounding each integral by
$C/(1+|t\mp n|)$ gives $\sum |a_n(t)|=O(\log(2+|t|))$ over
this range, and the weights $(1+|n|)^p$ are at most $(3+2|t|)^p$
there. For $|n|>2|t|+2$ direct integration gives
\[
 a_n(t)=\frac{t\bigl((-1)^ne^{i\pi t}-1\bigr)}{\pi i(t^2-n^2)}.
\]
Thus the weighted sum over $|n|>2|t|+2$ is bounded by
$C|t|\sum_{n>2|t|+2}n^{p-2}=O_p((1+|t|)^p)$. Translation preserves
the absolute values of the coefficients, proving the assertion for $c$.

Now take $u=\delta_g$ for $g\in S\cup S^{-1}$. For each $0<p<1$,
the bound on $\|u^n\|_{1,\omega}=\omega(g^n)$ gives a bounded unital
homomorphism $\pi_u:W_p\to\ell^1(G,\omega)$,
$f\mapsto\sum_n\widehat f(n)u^n$. The elements $h=\pi_u(b)$ and
$k=\pi_u(c)$ do not depend on $p$. For every $\varepsilon>0$, choose
$0<p<\min\{1,\varepsilon\}$. Applying $\pi_u$ to the estimated
exponentials gives
\[
 \|e^{ith}\|_{1,\omega}+\|e^{itk}\|_{1,\omega}
               =O_{g,\varepsilon}((1+|t|)^\varepsilon).
\]
For the ordinary group algebra one can take $p=0$ and obtain
$O(\log(2+|t|))$ directly. By \eqref{eq:circle-cosines},
$u=\cos h+i\cos k$. These pairs therefore generate the whole algebra:
the point masses of $S\cup S^{-1}$ generate the finitely supported
functions, which are dense. Corollary~\ref{cor:source-growth} applies.
\end{proof}

Section~\ref{sec:rank} compares this result for the free group on
countably many generators with the infinite exponential rank of its
group algebra.

\section{Factorisation rank and the number of tests}\label{sec:rank}

Let $G_0(A)$ be the identity component of the invertible group of $A$,
and let $E_m(A)$ consist of the products of $m$ exponentials.
Since a neighbourhood of $\one$ has a logarithm,
$G_0(A)=\bigcup_{m\ge1}E_m(A)$. In the terminology of Phillips and
Ringrose \cite[Definition~1.1]{PhRi}, $\operatorname{ber}(A)\le m$
means that $E_m(A)=G_0(A)$, whereas $\operatorname{ber}(A)\le m+\varepsilon$
means that $E_m(A)$ is norm-dense in $G_0(A)$. Here $m+\varepsilon$
is a formal value strictly between $m$ and $m+1$.
If $\operatorname{ber}(A)\le1$, then every product of two
exponentials is an exponential, so \eqref{hyp:E1} and \eqref{hyp:E2}
coincide and Theorem~\ref{thm:main} shows that a single exponential
suffices. Density alone does not give this, since a limit of invertible
elements need not be invertible.

For example, $\C\one+\mathcal K(X)$ has exponential rank one for every
infinite-dimensional complex Banach space $X$. This subalgebra is
inverse-closed in $\Bcal(X)$: an invertible $a=\lambda\one+K$ has
$\lambda\ne0$ and $a^{-1}=\lambda^{-1}(\one-a^{-1}K)$.
By Riesz--Schauder theory its spectrum is countable. If $a$ is
invertible, a ray from zero therefore misses its spectrum, so a
branch of the logarithm on the slit plane gives a logarithm of $a$
in the algebra.

The Hilbert space case illustrates the distinction between exact
factorisation and the number of tests. For infinite-dimensional complex $H$,
\[
                        \operatorname{ber}(\Bcal(H))=2 .
\]
Indeed, every invertible $a=u|a|$ is a product
$e^{ik}e^{\log|a|}$: the bounded Borel functional calculus gives a
bounded self-adjoint $k$ with $u=e^{ik}$.
For the lower bound, the Halmos--Lumer example gives a nonempty open
set of invertibles without square roots, hence without logarithms;
see \cite[Proposition~1.4(4) and its proof]{PhRi}.
The example on a separable summand $H_0$ also works on a larger
space: extend an interior point $v$ of that set by $cI$ on
$H_0^\perp$, where $c>\|v\|+1$. Under a small perturbation the Riesz
projection separating $\sigma(v)$ from $c$ remains close to the
projection onto $H_0$. Conjugating the two projections by an
invertible close to the identity identifies the perturbed restriction
with an operator close to $v$. A square root would commute with the
Riesz projection and restrict to a square root there. In finite
dimension, by contrast, every invertible matrix has a logarithm.

The polar decomposition gives the upper bound two for every von
Neumann algebra. The same bound passes to its quotients by proper
closed two-sided ideals:
a product of exponentials in the quotient lifts factor by factor to
the identity component upstairs. For the Calkin algebra on separable
infinite-dimensional Hilbert space, the exact rank
is again two \cite[Proposition~1.4(5)]{PhRi}.

The algebra $C(\T,M_2(\C))$ likewise has exponential rank two
\cite[p.~356]{PhillipsSurvey1994}. An example explaining the lower
bound is $v(\zeta)=\operatorname{diag}(2\zeta,\bar\zeta/2)$.
It belongs to $E_2$: writing
$J(s)=\bigl(\begin{smallmatrix}0&s\\s^{-1}&0\end{smallmatrix}\bigr)$,
we have $v(\zeta)=J(2\zeta)J(1)$, and every involution $J$ equals
$\exp(i\pi(\one-J)/2)$. For a small uniform perturbation of $v$, the
outer eigenvalue remains simple and close to $2\zeta$, so still has
winding number one. If the perturbed matrix function were $e^h$, its
rank-one outer Riesz projection $P$ would commute with $h$ and depend
continuously on $\zeta$. Then $\operatorname{tr}(hP)$, with the ordinary
matrix trace, would be a continuous logarithm of the outer eigenvalue,
which is impossible. Thus even a neighbourhood of $v$ misses $E_1$.

For discrete group algebras the two numbers can be as far apart as
possible.
For the free group $F_\infty$ on countably many generators, Phillips
proved in \cite[Corollary~7.5]{PhillipsSurvey1994} that
\[
                   \operatorname{ber}(\ell^1(F_\infty))=\infty .
\]
Indeed, the canonical homomorphism to the full group $C^*$-algebra
$C^*(F_\infty)$ has dense range. By approximating the individual
exponents, every finite exponential product downstairs is a limit
of images of such products upstairs. Thus density of $E_m$ in
$G_0(\ell^1(F_\infty))$ for a fixed $m$ would imply its density in
$G_0(C^*(F_\infty))$, contrary to \cite[Corollary~2.7]{PhillipsMany1994}.
Yet a single exponential suffices to test every unital linear
surjection from $\ell^1(F_\infty)$ onto a semisimple target, by
Corollary~\ref{cor:E1-group-source}. The number of factors needed to
represent invertibles can therefore be unbounded even when the
one-exponential preserver test suffices.

\section{What a single exponential detects in the target}\label{sec:targets}

\subsection{Tracial states and the square defect}
For a unital $C^*$-algebra $B$ put
\[
 I_{\rm tr}=\bigcap_{\tau}\{b\in B:\tau(b^*b)=0\},
\]
where $\tau$ runs over the tracial states; put $I_{\rm tr}=B$ if there
are none. Each tracial null space is a closed left ideal, by the
Cauchy--Schwarz inequality and $b^*c^*cb\le\|c\|^2b^*b$; it is
self-adjoint by traciality, and hence a two-sided ideal. Thus $I_{\rm tr}=0$
precisely when the tracial states form a faithful family.

\begin{theorem}\label{thm:tracial-E1}
Let $T:A\to B$ be a unital complex-linear surjection from a unital
Banach algebra onto a unital $C^*$-algebra. If $T(e^a)$ is invertible
for every $a\in A$, then
\[
            T(x^2)-(Tx)^2\in I_{\rm tr}\qquad(x\in A).
\]
In particular, a single exponential suffices when $B$ has a faithful
family of tracial states.
\end{theorem}

The proof uses the trace logarithm of Aupetit and Mathieu--Schulz
\cite{Aupetit1998Survey,MathieuSchulz2023} and a substitution used by
Christensen for matrix-valued maps
\cite[proof of Theorem~3.5]{Christensen1997}. It does not require the
spectral-isometry theorem.
\begin{proof}
Proposition~\ref{prop:E1-hull} gives continuity of $T$ and
$r_B(Tx)\le r_A(x)$. Fix a tracial state $\tau$ and set $f=\tau T$.
Since $|\tau(b)|\le r_B(b)$ \cite[Section~2]{MathieuSchulz2023},
the entire function $z\mapsto f(e^{za}xe^{-za})$ is bounded by
$r_A(x)$. Differentiating at zero gives $f(ax)=f(xa)$.

For $a\in A$, write $F(z)=T(e^{za})$ and let $g$ be the entire
primitive of $\tau(F'F^{-1})$ with $g(0)=0$. Then
\[
 \Re g(z)=\tfrac12\tau\bigl(\log(F(z)^*F(z))\bigr)
 \le\log\|F(z)\|\le\log\|T\|+|z|\,\|a\|.
\]
For the identity, along a smooth path of positive invertibles $P$ one
has $(\tau(\log P))'=\tau(P^{-1}P')$: differentiate the resolvent
integral for $\log P$, use cyclicity and integrate
$\int_0^\infty(P+s\one)^{-2}\,ds=P^{-1}$. Apply this to $P=F^*F$
along a path from zero. The inequality follows from positivity of
$\tau$. Borel--Carath\'eodory on discs of radii $R$ and $2R$, followed
by Cauchy's estimates, gives $g(z)=z g'(0)=z f(a)$.

Now let $a=\log(\one+t x)$ for $t$ sufficiently small. For
$0\le s\le1$, $F(s)$ is uniformly close to $\one$ and termwise
differentiation of the logarithm series gives
$\frac{d}{ds}\tau(\log F(s))=\tau(F'(s)F(s)^{-1})$. Integration in
$s$ therefore gives
\[
 \tau(\log(\one+tTx))=g(1)=f(\log(\one+t x)).
\]
Comparison of coefficients yields $\tau((Tx)^k)=f(x^k)$ for $k\ge1$.
Substitute $x+ty$ for $x$ when $k=2,3$ and use traciality on both
sides. The coefficients of $t$ give
\[
 f(xy)=\tau(Tx\,Ty),\qquad f(x^2y)=\tau((Tx)^2Ty).
\]
Replacing $x$ by $x^2$ in the first identity shows that
$\tau(d\,Ty)=0$ for $d=T(x^2)-(Tx)^2$ and every $y$. Choosing $y$ with
$Ty=d^*$, we obtain $\tau(dd^*)=0$, that is, $\tau(d^*d)=0$ by
traciality. As $\tau$ was arbitrary, $d\in I_{\rm tr}$.
\end{proof}

The theorem applies to finite von Neumann algebras. It also gives a
Banach algebra consequence, without the spectral-isometry theorem: if
$B$ has a separating family of finite-dimensional irreducible
representations, the single-exponential test forces $T$ to be Jordan.
Indeed, each such representation $\pi$ maps $B$ onto a full matrix
algebra, by Burnside's theorem, and $\pi T$ is a unital surjection onto
it satisfying \eqref{hyp:E1}. The normalised trace is faithful on a
matrix algebra, so Theorem~\ref{thm:tracial-E1} shows that $\pi T$ is
Jordan, and the representations separate the square defects of $T$.
For maps preserving all invertible elements this is Aupetit's theorem
\cite{AupetitPJM}.

\subsection{The primitive obstruction}
Call a linear map $S$ a \emph{vector-valued trace} if $S(ax)=S(xa)$
for all $a,x$; this condition makes no assumption on multiplication in
its range. A primitive Banach algebra has a faithful algebraically
irreducible representation and is prime and semisimple.

\begin{theorem}[The primitive alternative]\label{thm:primitive-E1}
Let $S:A\to P$ be a unital linear surjection onto a primitive Banach
algebra and suppose that $S(e^a)$ is invertible for all $a\in A$.
Then $S$ is a Jordan epimorphism or a vector-valued trace.
\end{theorem}
\begin{proof}
By Proposition~\ref{prop:E1-hull} and Lemma~\ref{lem:E1-conjugation},
$\ker S$ is a Lie ideal and $D_a(Sx)=S([a,x])$ defines a derivation
of $P$. Thus $\beta(Sa,Sx)=S([a,x])$ is well defined and alternating
(it is well defined in the first variable because
$S([a,x])=-S([x,a])$), and it is a derivation in each variable. If $\beta=0$, then $S$ is a
vector-valued trace. Otherwise the two Leibniz expansions of
$\beta(bd,ce)$ give
\[
 [b,c]\beta(d,e)=\beta(b,c)[d,e]\qquad(b,c,d,e\in P);
\]
compare \cite[equation~(5.2)]{AznayChabbabi2026}. Whenever
$\beta(b,c)=0$, this identity gives $[b,c]\beta(d,e)=0$ for all
$d,e$. Replacing $d$ by $wd$ and expanding $\beta(wd,e)$ once more, we
obtain $[b,c]w\beta(d,e)=0$ for all $w$. As some $\beta(d,e)\ne0$,
primeness gives $[b,c]=0$. In particular $S$ sends commuting pairs
to commuting pairs.

Given $x\in A^{-1}$, let $C$ be the closed unital subalgebra of $A$
generated by $x$ and $x^{-1}$; it is commutative. Its image $S(C)$
consists of pairwise commuting elements, and therefore lies in a
maximal commutative subalgebra $D$ of $P$. Such a $D$ is closed and inverse-closed: the
inverse of any of its elements invertible in $P$ still commutes with
$D$. Thus for each character $\chi$ of $D$, the functional
$\chi S|_C$ is unital and never vanishes on exponentials of $C$. By
\eqref{eq:scalar-exponential} it is multiplicative, so
$\chi(Sx)\chi(S(x^{-1}))=1$. No character of $D$ vanishes on $Sx$;
hence $Sx$ is invertible. The invertibility-preserver theorem
\cite[Proposition~5.2]{AznayChabbabi2026} now gives the Jordan identity.
\end{proof}

\begin{corollary}\label{cor:primitive-obstruction}
A counterexample to the single-exponential assertion exists if and
only if there is a unital surjective vector-valued trace satisfying
\eqref{hyp:E1} onto an infinite-dimensional noncommutative primitive
Banach algebra.
\end{corollary}
\begin{proof}
Let $T$ be a counterexample, that is, a unital linear surjection onto
a semisimple Banach algebra $B$ which satisfies \eqref{hyp:E1} and is
not a Jordan epimorphism. Some square defect of $T$ is nonzero, and it
survives in some primitive quotient $B/Q$, since the intersection of the
primitive ideals of $B$ is $\rad(B)=\{0\}$. The composition of $T$ with
the quotient map is a unital surjection onto $B/Q$ which satisfies
\eqref{hyp:E1} and is not Jordan, so it is a vector-valued trace by
Theorem~\ref{thm:primitive-E1}. A commutative primitive
target is $\C$, where \eqref{eq:scalar-exponential} rules out a defect, and
a finite-dimensional primitive target is a full matrix algebra, where
Theorem~\ref{thm:tracial-E1} rules it out; so neither can exhibit the
defect.

Conversely, such a trace $S$ is a counterexample unless it is Jordan;
suppose that it is. Polarisation gives
$S(xy)=S(x)\circ S(y)$, where $b\circ c=(bc+cb)/2$.
By surjectivity, $\circ$ is then associative on the range. Since
$(b\circ c)\circ d-b\circ(c\circ d)=[c,[b,d]]/4$, all commutators are
central. For $q=[b,c]$
the element $bq=[b,bc]$ is again a commutator, hence central; as $q$ is
central as well, $q^2=q(bc-cb)=bqc-cqb=[bq,c]=0$. A non-zero central element
whose square is zero generates a non-zero two-sided ideal with square zero,
so $q=0$ by semiprimeness. As $b$ and $c$ were arbitrary, the range would
then be commutative, which is impossible in a primitive algebra different
from $\C$.
\end{proof}

In particular, such a trace must factor as a bounded linear map through
the Banach-space quotient of $A$ by the closed span of its radical and
its commutators (Proposition~\ref{prop:E1-hull} and
Corollary~\ref{cor:E1-source-radical}). This quotient would have to be
infinite-dimensional.

The target $\Bcal(\ell^2)$ is an instructive case left open by these
results.
It is primitive and has no tracial states: the isometries onto the
even and odd coordinates have orthogonal range projections summing to
the identity, whereas a tracial state would give each projection
value one. Must a unital linear surjection from an arbitrary unital
Banach algebra onto $\Bcal(\ell^2)$ be Jordan if it sends single
exponentials to invertibles? By Theorem~\ref{thm:primitive-E1},
any counterexample onto this target would have to annihilate every
commutator.

\subsection*{Acknowledgements}
The author acknowledges the support of the CDP
C2EMPI and its institutional partners for the project R-CDP-24-004-C2EMPI,
and of the EU COST Actions.

\end{document}